\documentclass[12pt]{amsart}
\usepackage{amsfonts}
\usepackage{amssymb}
\usepackage{geometry}

\theoremstyle{definition}

\theoremstyle{remark}

\numberwithin{equation}{section}

\input{tcilatex}
\begin{document}
\title{Two Dimensional Representations of the Braid Group B3 and its Burau
Representation Squared}
\author{Mehmet K\i rdar}
\address{Department of Mathematics, Faculty of Arts and Science, Nam\i k
Kemal University, Tekirda\u{g}, Turkey}
\email{ mkirdar@nku.edu.tr}
\date{4 July, 2019 \textit{Mathematics Subject Classification. }[2010]
Primary 20F36; Secondary 20C99}
\keywords{Braid Group B3, Burau Representation}

\begin{abstract}
We list the irreducible two dimensional complex representations of the Braid
group B3 in elementary way. Then, we make a decomposition of the square of
its irreducible Burau representation.
\end{abstract}

\maketitle


\section{Introduction}

Artin braid group $B_{n}$, $n\geq 2$, is given by generators $\sigma _{1}$, $%
\sigma _{2}$, ..., $\sigma _{n-1}$ subject to the relations:

\begin{center}
$%
\begin{array}{l}
\left( 1\right) \text{ }\sigma _{i}\sigma _{j}=\sigma _{j}\sigma _{i}\text{
for }\left\vert i-j\right\vert \geq 2 \\ 
\left( 2\right) \text{ }\sigma _{i+1}\sigma _{i}\sigma _{i+1}=\sigma
_{i}\sigma _{i+1}\sigma _{i}\text{ for }1\leq i\leq n-2.%
\end{array}%
$
\end{center}

A quick proof of the fact that these are the minimal relations for the
geometric braid group is given by Jun Morita in [3] by induction.

Note that $B_{2}=Z$ is not interesting, on the other hand $B_{3}$ is
generated by $\sigma _{1}$, $\sigma _{2}$ subject to a single relation $%
\sigma _{1}\sigma _{2}\sigma _{1}=\sigma _{2}\sigma _{1}\sigma _{2}$ and it
is already very complicated. For $n\geq 3,$ Braid groups are non-abelian
infinite groups.

A $k$-dimensional complex representation of $B_{n}$ is a group homomorphism
from this group into $GL_{k}(%
\mathbb{C}
),$ the group of $k\times k$ complex invertible matrices. Two
representations are isomorphic if they are conjugate by a fixed invertible
matrix. We say that a representation is irreducible if it is not isomorphic
to direct sum of two lower dimensional representations. All 1-dimensional
representations are irreducible by definition. A famous example of a higher
dimensional irreducible representation of the braid group $B_{n}$ is its
irreducible Burau representation $\beta (z)$, [1].

Edward Formanek's paper [1] on Braid group representations of low degrees
discusses the classification of the complex representations of the Braid
group $B_{n}$ for dimensions $\leq n-1.$But, his first and easiest
classification result, Theorem 11, seems to list less 2-dimensional
representations of $B_{3}$ than we had by solving matrix equations.

In the first part of this note, we classify the complex 2-dimensional
representations of $B_{3}$ in the elementary way. Then, in the second part,
we will make the decomposition of $\beta (z)\otimes \beta (z)$ hoping to
find an analogous of the relation that we have in $R(S_{3}).$ It turned out
to be different, not 1+1+2 decomposition but 1+3. The 3-dimensioanal
irreducible representation we obtained obeys a Pascal triangle pattern. This
kind of patterns first observed in [2]. Higher dimensional versions should
be related to higher tensor powers of the Burau representation and shouldn't
be hard to generalize.

\section{Two Dimensional Representations of $B_{3}$}

Due to the relation $\sigma _{1}\sigma _{2}\sigma _{1}=\sigma _{2}\sigma
_{1}\sigma _{2}$, if $\xi :\sigma _{1}\rightarrow z_{1}$, $\sigma
_{2}\rightarrow z_{2}$ where $z_{1}$, $z_{2}\in 
\mathbb{C}
^{\ast }$ is a 1-dimensional representation of $B_{3}$, then we must have $%
z_{1}=z_{2}.$ So for each $z\in 
\mathbb{C}
^{\ast },$ let us denote the representation, $\sigma _{1},$ $\sigma
_{2}\rightarrow z$ by $\xi (z).$ All non-isomorphic 1-dimensional
representations of $B_{3}$ and also of other higher order braid groups are
in this form.

Next, we want to determine all possible forms of irreducible 2-dimensional
representations.

\bigskip

\textbf{Theorem 1.} All irreducible 2-dimensional representations of $B_{3}$
up to tensor with 1-dimansional representations are in the form

(i) $\sigma _{1}\rightarrow \left[ 
\begin{array}{cc}
-z & 0 \\ 
0 & 1%
\end{array}%
\right] ,\sigma _{2}\rightarrow \left[ 
\begin{array}{cc}
\frac{1}{z+1} & f \\ 
g & -\frac{z^{2}}{z+1}%
\end{array}%
\right] $ where $fg=\dfrac{z(z^{2}+z+1)}{(z+1)^{2}}$, $z\neq 0,-1$ and $%
z^{2}+z+1\neq 0,$

or in the form

(ii) $\sigma _{1}\rightarrow \left[ 
\begin{array}{cc}
1 & z \\ 
0 & 1%
\end{array}%
\right] ,\sigma _{2}\rightarrow \left[ 
\begin{array}{cc}
e & z(e-1)^{2} \\ 
-\frac{1}{z} & 2-e%
\end{array}%
\right] $ where $z\neq 0.$

\textit{Proof.}\textbf{\ }Let us now consider a 2-dimensional representation
of $B_{3}$ given by $\sigma _{1}\rightarrow A$, $\sigma _{2}\rightarrow B$
where $A$, $B\in GL_{2}(%
\mathbb{C}
).$

Since we are interested in isomorphism classes of 2-dimensional
representations, we can diagonalize $A$ if it is diagonalizable or if not,
we can put it in Jordan form. So, we have two cases.

(i) We assume that $A$ is diagonal and also we can tensor with a
1-dimensional representation to assume that one entry of the diagonal is 1.

So, let us assume that $A=\left[ 
\begin{array}{cc}
-z & 0 \\ 
0 & 1%
\end{array}%
\right] $ \ and that $B=\left[ 
\begin{array}{cc}
e & f \\ 
g & h%
\end{array}%
\right] $ where $z\neq 0$. The choice of $-z$ in the matrix $A$ will be
clarified in the next section.

Due to the relation $\sigma _{1}\sigma _{2}\sigma _{1}=\sigma _{2}\sigma
_{1}\sigma _{2},$ we must have $ABA=BAB$ and thus we can solve for $e,h$ and
the product $fg$ to get the result as expressed above.

Since if $z^{2}+z+1\neq 0$ the given matrices do not have common
eigenvectors for their same eigenvalues, we can not find one dimensional
subrepresentation. So, the representations given in this case are all
irreducible.

(ii) We assume that $A$ is in Jordan from and also we can tensor with a
1-dimensional representation to assume that diagonal entries are 1.

So let $A=\left[ 
\begin{array}{cc}
1 & z \\ 
0 & 1%
\end{array}%
\right] $ \ and that $B=\left[ 
\begin{array}{cc}
e & f \\ 
g & h%
\end{array}%
\right] $ where $z\neq 0$. Due to the relation $\sigma _{1}\sigma _{2}\sigma
_{1}=\sigma _{2}\sigma _{1}\sigma _{2},$ we must have $ABA=BAB$ and thus we
find that $f=z(e-1)^{2},$ $g=-\dfrac{1}{z}$ and $h=2-e$.

Since we have one eigenvector in this case, we can not split these
representations. So, all representations in this case are irreducible too.

Edward Formanek's paper [1] also discusses the condition $z^{2}+z+1\neq 0$
related to Burau representation. He generalizes this condition for all $n$
in Lemma 6 of the paper. His main tool is pseudoreflections.

\section{Burau Representation of $B_{3}$ Squared}

The irreducible Burau representation of $B_{3}$ is given by

\begin{equation*}
\sigma _{1}\rightarrow \left[ 
\begin{array}{cc}
-z & 0 \\ 
1 & 1%
\end{array}%
\right] ,\sigma _{2}\rightarrow \left[ 
\begin{array}{cc}
1 & z \\ 
0 & -z%
\end{array}%
\right]
\end{equation*}%
where $z\in 
\mathbb{C}
^{\ast }$ is a parameter. We will denote this representation by $\beta (z).$
Hence, $\beta (z):B_{3}\rightarrow GL_{2}(%
\mathbb{C}
).$

\bigskip

\textbf{Proposition 2. }The irreducible Burau representation of $B_{3}$ is
isomorphic to the representation, $z\neq 1,$

\begin{equation*}
\sigma _{1}\rightarrow \left[ 
\begin{array}{cc}
-z & 0 \\ 
0 & 1%
\end{array}%
\right] ,\sigma _{2}\rightarrow \left[ 
\begin{array}{cc}
\frac{1}{z+1} & -\frac{z}{z+1} \\ 
-\frac{z^{2}+z+1}{z+1} & -\frac{z^{2}}{z+1}%
\end{array}%
\right]
\end{equation*}

\textit{Proof. }The matrix $A=\left[ 
\begin{array}{cc}
-z & 0 \\ 
1 & 1%
\end{array}%
\right] $ is diagonalizable by $P=\left[ 
\begin{array}{cc}
-(z+1) & 0 \\ 
1 & 1%
\end{array}%
\right] $ to $P^{-1}AP=\left[ 
\begin{array}{cc}
-z & 0 \\ 
0 & 1%
\end{array}%
\right] .$ And then for $B=\left[ 
\begin{array}{cc}
1 & z \\ 
0 & -z%
\end{array}%
\right] ,$ we get $P^{-1}BP=\left[ 
\begin{array}{cc}
\frac{1}{z+1} & -\frac{z}{z+1} \\ 
-\frac{z^{2}+z+1}{z+1} & -\frac{z^{2}}{z+1}%
\end{array}%
\right] $ as required.

We notice that the Burau representation of $B_{3}$ can be obtained by taking 
$f=-\frac{z}{z+1}$ in Theorem 1 (i).

$S_{3},$ the symmetric group on 3 letters is generated by cycles $s_{1}$ and 
$s_{2}$ subject to the relations $s_{1}^{2}=s_{2}^{2}=1$ and $%
s_{1}s_{2}s_{1}=s_{2}s_{1}s_{2}.$

We have the natural projection homomorphism $\pi :B_{3}\rightarrow S_{3}$
which sends $\sigma _{1}$ to $s_{1}$ and $\sigma _{2}$ to $s_{2}$. Note that
we have these natural homomorphisms $\pi :B_{n}\rightarrow S_{n}$ for each $%
n $ and the kernel of these homomorphisms are called pure braid groups.

The complex representation ring of $S_{3}$ is given by

\begin{equation*}
R(S_{3})=%
\mathbb{Z}
\lbrack \xi ,\rho ]/(\xi ^{2}=1,\text{ }\rho ^{2}=\rho +\xi +1,\text{ }\xi
\rho =\rho )
\end{equation*}

where $\xi $ is the 1-dimensional representation given by $s_{1}\rightarrow
-1,$ $s_{2}\rightarrow -1$ and is called the sign representation and $\rho $
is the 2-dimensional representation given by $s_{1}\rightarrow \left[ 
\begin{array}{cc}
-1 & 0 \\ 
1 & 1%
\end{array}%
\right] ,s_{2}\rightarrow \left[ 
\begin{array}{cc}
1 & 1 \\ 
0 & -1%
\end{array}%
\right] $ which is called the standard representation.

The standard representation of $S_{3}$ is obtained from Burau representation
of $B_{3}$by taking $z=1.$ In other words, $\rho =\pi \circ \beta (1).$
Inspired of this connection, we suspect that $\beta (z)\otimes \beta (z)$
might satisfy a similar relation.

\bigskip

\textbf{Proposition 3. }Let $z\neq 0,-1.$\textbf{\ }Then $\beta (z)^{2}=\xi
(-z)+\mu (z)$ where $\mu (z)$ is a 3-dimensional representation which is
irreducible if $z^{3}\neq 1$ and it is given by

$\sigma _{1}\rightarrow \left[ 
\begin{array}{ccc}
1 & 0 & 0 \\ 
0 & -z & 0 \\ 
0 & 0 & z^{2}%
\end{array}%
\right] ,$ $\sigma _{2}\rightarrow \left[ 
\begin{array}{ccc}
\frac{z^{4}}{\left( z+1\right) ^{2}} & \frac{z^{2}}{\left( z+1\right) ^{2}}%
\left( z^{2}+z+1\right) & \frac{1}{\left( z+1\right) ^{2}}\left(
z^{2}+z+1\right) ^{2} \\ 
2\frac{z^{3}}{\left( z+1\right) ^{2}} & z\frac{z^{2}+1}{\left( z+1\right)
^{2}} & -\frac{2}{\left( z+1\right) ^{2}}\left( z^{2}+z+1\right) \\ 
\frac{z^{2}}{\left( z+1\right) ^{2}} & -\frac{z}{\left( z+1\right) ^{2}} & 
\frac{1}{\left( z+1\right) ^{2}}%
\end{array}%
\right] $

\textit{Proof. }First, we calculate

$A=\beta (z)\otimes \beta (z)(\sigma _{1})=$ $\left[ 
\begin{array}{cc}
-z & 0 \\ 
1 & 1%
\end{array}%
\right] \otimes \left[ 
\begin{array}{cc}
-z & 0 \\ 
1 & 1%
\end{array}%
\right] =\left[ 
\begin{array}{cccc}
z^{2} & 0 & 0 & 0 \\ 
-z & -z & 0 & 0 \\ 
-z & 0 & -z & 0 \\ 
1 & 1 & 1 & 1%
\end{array}%
\right] $ and that

$B=\beta (z)\otimes \beta (z)(\sigma _{2})=$ $\left[ 
\begin{array}{cc}
1 & z \\ 
0 & -z%
\end{array}%
\right] \otimes \left[ 
\begin{array}{cc}
1 & z \\ 
0 & -z%
\end{array}%
\right] =\left[ 
\begin{array}{cccc}
1 & z & z & z^{2} \\ 
0 & -z & 0 & -z^{2} \\ 
0 & 0 & -z & -z^{2} \\ 
0 & 0 & 0 & z^{2}%
\end{array}%
\right] .$

Now, $A$ is diagonalizable with the matrix $P=\left[ 
\begin{array}{cccc}
0 & 0 & 0 & z^{2}+2z+1 \\ 
0 & -z-1 & -1 & -z-1 \\ 
0 & 0 & 1 & -z-1 \\ 
1 & 1 & 0 & 1%
\end{array}%
\right] $ and

$P^{-1}AP=$ $\left[ 
\begin{array}{cccc}
1 & 0 & 0 & 0 \\ 
0 & -z & 0 & 0 \\ 
0 & 0 & -z & 0 \\ 
0 & 0 & 0 & z^{2}%
\end{array}%
\right] .$

And then $P^{-1}BP=\allowbreak \allowbreak \left[ 
\begin{array}{cccc}
\frac{z^{4}}{\left( z+1\right) ^{2}} & \frac{z^{2}}{\left( z+1\right) ^{2}}%
\left( z^{2}+z+1\right) & 0 & \frac{1}{\left( z+1\right) ^{2}}\left(
z^{2}+z+1\right) ^{2} \\ 
2\frac{z^{3}}{\left( z+1\right) ^{2}} & z\frac{z^{2}+1}{\left( z+1\right)
^{2}} & 0 & -\frac{2}{\left( z+1\right) ^{2}}\left( z^{2}+z+1\right) \\ 
-\frac{z^{3}}{z+1} & -\frac{z}{z+1}\left( z^{2}+z+1\right) & -z & \frac{1}{%
z+1}\left( z^{2}+z+1\right) \\ 
\frac{z^{2}}{\left( z+1\right) ^{2}} & -\frac{z}{\left( z+1\right) ^{2}} & 0
& \frac{1}{\left( z+1\right) ^{2}}%
\end{array}%
\right] \allowbreak .$

Let $%
\mathbb{C}
^{4}=Span\left\{ e_{1},e_{2},e_{3},e_{4}\right\} .$ From the calculations
above, we observe that $Span\left\{ e_{3}\right\} $ is an invariant subspace
of the representation $\beta (z)\otimes \beta (z)$ with eigenvalue $-z.$ So $%
\beta (z)\otimes \beta (z)$ is reducible. Hence, we have $\beta (z)\otimes
\beta (z)=\xi (-z)\oplus \mu (z)$ where $\xi (-z)$ is the mentioned
1-dimensional subrepresentation and

$\mu (z):$

$\sigma _{1}\rightarrow C=\left[ 
\begin{array}{ccc}
1 & 0 & 0 \\ 
0 & -z & 0 \\ 
0 & 0 & z^{2}%
\end{array}%
\right] ,$

$\sigma _{2}\rightarrow D=\left[ 
\begin{array}{ccc}
\frac{z^{4}}{\left( z+1\right) ^{2}} & \frac{z^{2}}{\left( z+1\right) ^{2}}%
\left( z^{2}+z+1\right) & \frac{1}{\left( z+1\right) ^{2}}\left(
z^{2}+z+1\right) ^{2} \\ 
2\frac{z^{3}}{\left( z+1\right) ^{2}} & z\frac{z^{2}+1}{\left( z+1\right)
^{2}} & -\frac{2}{\left( z+1\right) ^{2}}\left( z^{2}+z+1\right) \\ 
\frac{z^{2}}{\left( z+1\right) ^{2}} & -\frac{z}{\left( z+1\right) ^{2}} & 
\frac{1}{\left( z+1\right) ^{2}}%
\end{array}%
\right] $

is a 3-dimensional subrepresentation.$\allowbreak $

The matrix $D$ has eigenvectors $\left[ 
\begin{array}{c}
1 \\ 
-2 \\ 
1%
\end{array}%
\right] $ for the eigenvalue $1,$ $\left[ 
\begin{array}{c}
-\frac{1}{z}\left( z^{2}+z+1\right) \\ 
\frac{1}{z}\left( z^{2}+1\right) \\ 
1%
\end{array}%
\right] $ for the eigenvalue $-z,\allowbreak \left[ 
\begin{array}{c}
\frac{1}{z^{2}}\left( z^{4}+2z^{3}+3z^{2}+2z+1\right) \\ 
\frac{1}{z}\left( 2z^{2}+2z+2\right) \\ 
1%
\end{array}%
\right] $ for the eigenvalue $z^{2}\allowbreak .$ Whereas, eigenvectors of $%
C $ are the standard basis vectors of $%
\mathbb{C}
^{3}$. If $z\neq 1$ and $z^{2}+z+1\neq 0,$ which means that $z^{3}\neq 1,$
there is no common eigenvectors, and we deduce that these eigenvectors do
not produce a 1-dimensional subrepresentation. Therefore, $\mu (z)$ is
irreducible.

\bigskip

\textbf{Remark 4.} For $z=1,$ the eigenvectors of $D$ corresponding to 1 are 
$\left[ 
\begin{array}{c}
1 \\ 
-2 \\ 
1%
\end{array}%
\right] $ and $\left[ 
\begin{array}{c}
9 \\ 
6 \\ 
1%
\end{array}%
\right] $ and thus $\left[ 
\begin{array}{c}
3 \\ 
0 \\ 
1%
\end{array}%
\right] $ also is an eigenvector of $D$. $\left[ 
\begin{array}{c}
3 \\ 
0 \\ 
1%
\end{array}%
\right] $ is also an eigenvector of $C$ for $z=1$ corresponding to the
eigenvalue 1. This vector spans the trivial 1-dimensional subrepresentation $%
\xi (1)$ which corresponds to 1 in the decomposition $\rho ^{2}=\rho +\xi +1$
which occurred in $R(S_{3})$.

\bigskip

\textbf{Remark 5. }If we choose basis vectors, $\left\{ u_{1}=e_{1},\text{ }%
u_{2}=e_{2}+e_{3},\text{ }u_{3}=e_{4}\right\} ,$ we can easily show that the
representation $\mu (z)$ is isomorphic to the representation $\sigma
_{1}\rightarrow \left[ 
\begin{array}{ccc}
z^{2} & 0 & 0 \\ 
-z & -z & 0 \\ 
1 & 2 & 1%
\end{array}%
\right] ,$

$\sigma _{2}\rightarrow D=\left[ 
\begin{array}{ccc}
1 & 2z & z^{2} \\ 
0 & -z & -z^{2} \\ 
0 & 0 & z^{2}%
\end{array}%
\right] .$ We observe a Pascal triangle pattern in matrices. These patterens
also mentioned in [2].

$\allowbreak $

\end{document}